\documentclass[12pt]{article}
\usepackage{amsxtra,amssymb,amsthm,amsmath,latexsym}

\theoremstyle{plain}

\newtheorem{theorem}{Theorem}[section]

\newtheorem*{claim}{Claim}

\newtheorem{remark}[theorem]{Remark}
\numberwithin{equation}{section}
\newtheorem{remark2.1}{Remark 2.1}

\def\R{{\mathbb R}}

\def\oH{\buildrel\circ\over H}
\def\oH1{\buildrel\circ\over H\kern-.02in{}^1}

\def\supp{\hbox{\,supp\,}}
\def\const{\hbox{\,const\,}}
\def\ord{\hbox{\,ord\,}}

\begin{document}


\title{
Analytical solution of a new class of integral equations
}

\author{
A.G. Ramm\\
LMA/CNRS, 31 chemin Joseph Aiguier,
Marseille 13402, cedex 20, France\\
and mathematics department, KSU, Manhattan, KS 66506, USA\\
ramm@math.ksu.edu\\
http://www.math.ksu.edu/\,$\widetilde{\ }$\,ramm
}

\date{}

\maketitle\thispagestyle{empty}

\begin{abstract}
\footnote{Math subject classification:  45A05, 45H05, 93E10, 93E11, 34A30}
\\

Let $(1) Rh=f$, $0\leq x\leq L$, $Rh=\int^L_0 R(x,y)h(y)\,dy$,
where the kernel $R(x,y)$ satisfies the equation
$QR=P\delta(x-y)$. Here $Q$ and $P$ are formal differential
operators of order $n$ and $m<n$, respectively, $n$ and $m$
are nonnegative even integers, $n>0$, $m\geq 0$,
$Qu:=q_n(x)u^{(n)} + \sum^{n-1}_{j=0} q_j(x) u^{(j)}$,
$Ph:=h^{(m)} +\sum^{m-1}_{j=0} p_j(x) h^{(j)}$,
$q_n(x)\geq c>0$,
the coefficients $q_j(x)$ and $p_j(x)$ are smooth functions defined on
$\R$, $\delta(x)$ is the delta-function, $f\in H^\alpha(0,L)$,
$\alpha:=\frac{n-m}{2}$, $H^\alpha$ is the Sobolev space.

An algorithm for finding analytically the unique solution
$h\in\dot H^{-\alpha} (0,L)$ to (1) of minimal order of singularity is 
given.
Here $\dot H^{-\alpha}(0,L)$ is the dual space to $H^\alpha(0,L)$
with respect to the inner product of $L^2(0,L)$.

Under suitable assumptions it is proved that $R:\dot H^{-\alpha}(0,L) \to 
H^\alpha(0,L)$
is an isomorphism.

Equation (1) is the basic equation of random processes estimation theory.
Some of the results are generalized to the case of multidimensional
equation (1), in which case this is the basic equation of random
fields estimation theory.
\end{abstract}


\section{Introduction}
In monograph \cite{R1} (see also \cite{R3}) estimation theory for random 
fields and processes
is constructed. The estimation problem for a random process
is as follows. Let $u(x)=s(x)+n(x)$ be a random process observed on
the interval $(0,L)$, $s(x)$ is a useful signal and $n(x)$ is noise.
Without loss of generality we assume that $\overline{s(x)} = 
\overline{n(x)}=0$,
where the overbar stands for the mean value, 
$\overline{u^\ast(x)u(y)}:=R(x,y)$, $R(x,y)=\overline{R(y,x)}$,
$\overline{u^\ast(x)s(y)}:=f(x,y)$, and the star here stands for complex 
conjugate. 
The covariance functions $R(x,y)$ and $f(x,y)$ are assumed known.
One wants to estimate $s(x)$ optimally in the sense of minimum of the
variance of the estimation error. More precisely, one seeks a
linear estimate
\begin{equation}
Tu=\int^L_0 h(x,y) u(y)\,dy,
\end{equation}
such that
\begin{equation}
\overline{\left| (Tu)(x)-s(x)\right|^2} =\hbox{\ min.}
\end{equation}
This is a filtering problem. Similarly one can formulate the problem
of optimal estimation of $(As)(x)$, where $A$ is a known operator acting 
on $s(x)$. If $A=I$, where $I$ is the identity operator, then one has
the filtering problem, if $A$ is the differentiation operator, then
one has the problem of optimal estimation of the derivative of $s$,
if $As=s(x+x_0)$, then one has an extrapolation problem, etc.
The kernel $h(x,y)$ 
is, in general, a distribution. As in \cite{R1}, one derives a
necessary condition for $h$ to satisfy (1.2):
\begin{equation}
\int^L_0 R(x,y) h(y,z)\,dy=f(x,z), \quad 0\leq x,z\leq L.
\end{equation}
Since $z$ enters as a parameter in (1.3), the basic equation of
estimation theory is:
\begin{equation}
Rh:=\int^L_0 R(x,y)h(y)\,dy=f(x), \quad 0\leq x\leq L.
\end{equation}
The operator in $L^2(0,L)$ defined by (1.4) is symmetric.
In \cite{R1} it is assumed that the kernel
\begin{equation}
R(x,y)=\int^\infty_{-\infty} \frac{P(\lambda)}{Q(\lambda)} \Phi
(x,y,\lambda)\,d\rho(\lambda),
\end{equation}
where $P(\lambda)$ and $Q(\lambda)$ are positive polynomials,
$\Phi(x,y,\lambda)$ and $d\rho(\lambda)$ are spectral kernel and,
respectively, spectral measure of a selfadjoint
ordinary differential operator $\ell$ in $L^2(\R)$,
$\deg Q(\lambda)=q$, $\deg P(\lambda)=p<n$, $p\geq 0$,
$\ord \ell:=\sigma>0$. Actually in  \cite{R1} the multidimensional
case, when $x,y\in \R^r$, $r>1$, and $\ell$ is a selfadjoint elliptic 
operator in $L^2(\R^r)$ is also considered. 

It is proved in \cite{R1} that the operator
$R:\dot H^{-\alpha}(0,L) \to H^\alpha(0,L)$,
$\alpha:=\frac{n-m}{2}\sigma$
is an isomorphism. By $H^\alpha(0,L)$ the Sobolev space $W^{\alpha,2}(0,L)$ 
is denoted, and $\dot H^{-\alpha}(0,L)$ is the dual space to $H^\alpha(0,L)$
with respect to $L^2(0,L):=H^0(0,L)$ inner product. Namely,
$\dot H^{-\alpha}(0,L)$ is the space of distributions $h$ which are linear
bounded functionals on $H^\alpha(0,L)$.
The norm of $h\in\dot H^{-\alpha}(0,L)$ is given by the formula
\begin{equation}
\|h\|_{\dot H^{-\alpha}(0,L)}=\sup_{g\in H^\alpha(0,L)}
\frac{|(h,g)|}{\|g\|_{H^\alpha(0,L)}},
\end{equation}
where $(h,g)$ is the $L^2(0,L)$ inner product if $h$ and $g$
belong to $L^2(0,L)$.

One can also define  $\dot H^{-\alpha}(0,L)$ as
the subset of the elements of $H^{-\alpha}(\R)$ with support in $[0,L]$.

{\it In this paper we generalize the class of kernels $R(x,y)$
defined in (1.5):  we do not use the spectral theory, do not assume
$\ell$ to be selfadjoint, and do not assume that the operators $Q$ and $P$
commute.}

We assume that
\begin{equation}
QR=P\delta(x-y),
\end{equation}
where $Q$ and $P$ are formal differential operators of orders $n$ and $m$
respectively, $n>m\geq 0$, $n$ and $m$ are even integers, $\delta(x)$
is the delta-function, 
%
\begin{equation}
Qu:=\sum^n_{j=0}q_j(x)u^{(j)}, \quad q_n(x)\geq c>0, \quad Ph:=h^{(m)} 
+\sum^{m-1}_{j=0} p_j(x)h^{(j)},
\end{equation}
where $q_j$ and $p_j$ are smooth functions defined on $\R$.
We also assume that the equation $Qu=0$ has $ \frac{n}{2}$
linearly independent solutions $u^-_j \in L^2(-\infty,0)$ and 
$\frac{n}{2}$
linearly independent solutions $u^+_j \in L^2(0,\infty)$.
In particular, this implies that if $Qh=0$, $h\in H^\alpha(\R)$,
$\alpha>0$, then $h=0$, and the same conclusion holds 
for  $h\in H^\beta(\R)$ for any fixed real number $\beta$, including
negative  $\beta$, because any solution to the equation 
$Qh=0$ is smooth: it is a linear combination of
$n$ linearly independent solution to this equation, each of which is 
smooth and none belongs to $L^2(\R)$.

Let us assume that $R(x,y)$ is a selfadjoint kernel such that
\begin{equation}
c_1 \|\varphi\|^2_- \leq (R\varphi, \varphi)\leq c_2 \|\varphi\|^2_-,
\quad c_1=\const>0, \quad \forall \varphi\in C^\infty_0(\R),
\end{equation}
where $(\cdot,\cdot)$ is the $L^2(\R)$ inner product,
$\|\varphi\|_-:=\|\varphi\|_{H^{-\alpha}(\R)}:=\|\varphi\|_{-\alpha}$,
$\alpha:=\frac{n-m}{2}$, 
$\|\varphi\|_{\beta}:=\|\varphi\|_{H^{\beta}(\R)},$
and we use below the notation
$\|\varphi\|_+:=\|\varphi\|_{H^\alpha(0,L)}:=\|\varphi\|_{H^+}$.
The spaces $H^\alpha (0,L)$ and $\dot H^{-\alpha} (0,L)$ are dual of each 
other with respect to the $L^2(0,L)$ inner product, as was mentioned 
above.
If $\varphi\in \dot H^{-\alpha}(0,L)$, then
$\varphi\in H^{-\alpha}(\R),$ and the inequality (1.9) holds
for such $\varphi$. By this reason we also use
(for example, in the proof of Theorem 1.1 in Section 2) the notation
$H^-$ for the space $\dot H^{-\alpha} (0,L)$.

Assumption (1.9) holds, for example, if the following inequalities
(1.10) and (1.11) hold:
\begin{equation}
c_3\|\varphi\|_{-\alpha+n}
  \leq \|Q^\ast \varphi\|_{-\alpha}
  \leq c_4\|\varphi\|_{-\alpha+n},
  \ c_3,\ c_4=\const>0,
  \ \forall \varphi\in C^\infty_0(\R),
\end{equation}
\begin{equation}
c_5\|\varphi\|^2_{\frac{n+m}{2}} \leq (PQ^\ast \varphi,\varphi)
  \leq c_6\|\varphi\|^2_{\frac{n+m}{2}},
  \quad \forall \varphi\in C^\infty_0(\R),
\end{equation}
where $Q^\ast$ is a formally adjoint to $Q$ differential expression,
and $c_5$ and $c_6$ are positive constants independent of $\varphi\in 
C^\infty_0(\R)$. 
The right inequality (1.11) is obvious because $\ord PQ^\ast=n+m$, and
the right inequality (1.10) is obvious because $\ord Q^\ast=n.$

Let us formulate basic results of this paper.

\begin{theorem}
If (1.9) holds, then the operator $R$, defined in (1.5), is an isomorphism
of $\dot H^{-\alpha}(0,L)$ onto $H^\alpha(0,L)$, $\alpha=\frac{n-m}{2}$.
\end{theorem}

\begin{theorem}
If (1.7), (1.10) and (1.11) hold, then (1.9) holds and
$R:\dot H^{-\alpha}(0,L)\to H^\alpha(0,L)$ is an isomorphism.
\end{theorem}

\begin{theorem}
If (1.7), (1.10) and (1.11) hold,
and $f\in H^n(0,L)$, then
the solution to (1.4) in $\dot H^{-\alpha}(0,L)$ does exist,
is unique, and can be calculated analytically by the following formula:
\begin{equation}
h=\int^x_0 G(x,y)Q f\,dy + \sum^{n-\alpha-1}_{j=0}
\left[ a^-_j(-1)^j G^{(j)}_y (x,0) + a^+_j(-1)^j G^{(j)}_y(x,L) \right],
\end{equation}
where $a^\pm_j$ are some constants and $G(x,y)$ is the unique solution to
the problem
\begin{equation}
PG=\delta(x-y), \quad G(x,y)=0\ \hbox{for}\ x<y.
\end{equation}
The constants $a^\pm_j$ are uniquely determined from the condition
$h(x)=0$ for $x>L$.
\end{theorem}

\begin{remark}
The solution $h\in\dot H^{-\alpha}(0,L)$ is the solution to equation (1.4)
of  minimal order of singularity  (see \cite{R1}).
\end{remark}

\begin{remark}
If $P=1$ in (1.7) then the solution $h$ to (1.4) of minimal order of
singularity, $h\in\dot H^{- \frac{n}{2}}(0,L),$ can be calculated by the
formula $h=QF$, where $F$ is given by (2.9) (see below) and
$u_+$ and $u_-$ are the unique solutions of the problems $Qu_+=0$ if
$x>L$, $u^{(j)}_+(L)=f^{(j)}(L)$, $0\leq j\leq  \frac{n}{2}-1$,
$u_+(\infty)=0$, and $Qu_-=0$ if $x<0$, $u^{(j)}_-(0)=f^{(j)}(0)$,
$0\leq j\leq  \frac{n}{2}-1$, $u_-(-\infty)=0$.
\end{remark}

\section{Proofs}

\begin{proof}[Proof of Theorem 1.1.]
The set $C^\infty_0(0,L)$ is dense in $\dot H^{-\alpha} (0,L)$
(in the norm of  $ H^{-\alpha}(\R)$).
Using 
the right inequality (1.9),
one gets: 
%
\begin{equation}
\|R\|_{H^-\to H^+} = \sup_{h\in H^-} \frac{(Rh,h)}{\|h\|^2_-} \leq c_2,
\end{equation}
by the symmetry of $R$ in $L^2(0,L)$ (cf \cite[p.40]{R1}). This implies 
$||R||_{H^-\to H^+}\leq c_2$.
Using the left inequality (1.9), one gets:
$c_1\|h\|^2_- \leq \|Rh\|_+ \|h\|_-$, so
\begin{equation}
c_1\|h\|_- \leq \|Rh\|_+.
\end{equation}
Therefore
\begin{equation}
\|R^{-1}\|_{H^+\to H^-} \leq \frac{1}{c_1}.
\end{equation}
Consequently, the range $Ran(R)$ of $R$ is a closed subspace of $H^+.$ In 
fact,
$Ran(R)=H^+$. Indeed, if $Ran(R)\not= H^+$, then there exists a $g\in H^-$ 
such
that $0=(R\psi,g)$ $\forall \psi\in H^-$.
Taking $\psi=g$ and using the left inequality (1.9) one gets $\|g\|_-=0$, so
$g=0$. Thus $Ran(R)=H^+$.

Theorem 1.1 is proved.
\end{proof}

\begin{proof}[Proof of Theorem 1.2.]
From (1.7) and (1.8) it follows that the kernel $R(x,y)$ defines
a pseudodifferential operator of order $-2\alpha=m-n$.
In particular, this implies the right inequality (1.9).
In this argument inequalities (1.10) and (1.11) were not used.

Let us prove that (1.10) and (1.11) imply the left inequality (1.9).

One has
\begin{equation}
\|Q^\ast\varphi\|_{-\alpha} \leq C\|\varphi\|_{n-\alpha},
\quad \forall \varphi \in C^\infty_0(\R),
\end{equation}
because $\ord Q^\ast=n$. Inequality (1.10)  reads:
\begin{equation}
c_3\|\varphi\|_{-\alpha+n}\leq ||Q^\ast\varphi||_{-\alpha} \leq 
c_4\|\varphi\|_{-\alpha+n},
\quad \forall \varphi\in C^\infty_0(\R),
\end{equation}
where $c_3$ and $c_4$ are positive constants. If (2.5) holds,
then $Q^\ast: H^{-\alpha+n}(\R)\to H^{-\alpha}(\R)$
is an isomorphism of $H^{-\alpha+n}(\R)$ onto $H^{-\alpha}(\R)$
provided that $N(Q):=\{ w: Qw=0,\ w\in H^\alpha(\R)\} = \{0\}$.
Indeed, if the range of $Q^\ast$ is not all of $H^{-\alpha}(\R)$,
then there exists an $w\not= 0$, $w\in H^\alpha(\R)$
such that $(Q^\ast \varphi,w)=0$ \ $\forall \varphi\in C^\infty_0(\R)$,
so $Qw=0.$ If  $Qw=0$ and $w\in H^\alpha(\R)$, then, as was mentioned 
below formula (1.8), it follows that
 $w=0$. This proves that 
$Ran(Q^\ast)=H^{-\alpha}(\R)$.

Inequality (1.11) is necessary for the left inequality (1.9)
to hold. Indeed, let $\psi=Q^\ast\varphi$, $\varphi\in C^\infty_0(\R)$,
then (1.9) implies
\begin{equation}
c_5\|\varphi\|^2_{-\alpha+n} \leq c \|Q^\ast\varphi\|^2_{-\alpha}
\leq (RQ^\ast \varphi, Q^\ast \varphi)=(QRQ^\ast \varphi,\varphi) = 
(PQ^\ast\varphi,\varphi),
\end{equation}
where $c>0$ here (and elsewhere in this paper) stands for various 
estimation constants.
Because $-\alpha+n=\frac{n+m}{2}$, inequality (2.6) is the left
inequality (1.11). The right inequality (1.11) is obvious because
the order of the operator $PQ^\ast$ equals to $n+m$.

Let us prove now that inequalities (1.11) and (1.10) are sufficient for 
the left inequality (1.9) to hold.

Using the right inequality (1.10) and
the left inequality (1.11), one gets:
\begin{equation}
c\|\psi\|^2_{-\alpha} \leq c_5\|\varphi\|^2_{\frac{n+m}{2}}
\leq (PQ^\ast\varphi, \varphi)=(R\psi,\psi),
\quad \psi=Q^\ast\varphi,
\quad \forall \varphi\in C^\infty_0(\R).
\end{equation}
Let us prove that the set
$\{\psi\}=\{Q^\ast\varphi\}_{\forall \varphi\in C^\infty_0(\R)}$
is dense in $\dot H^{-\alpha}(0,L)$.
%
%
Assume the contrary. Then there is an $h\in \dot H^{-\alpha}(0,L)$,
$h\neq 0$,
such that $(Q^\ast\varphi, h)=0$ for all $\varphi\in C^\infty_0(\R)$.
Thus, $(\varphi, Qh)=0$ for all $\varphi\in C^\infty_0(\R)$. 
Therefore $Qh=0$, and, by the argument given below formula (1.8),
it follows that $h=0$.
This contradiction proves that the set 
$\{Q^\ast\varphi\}_{\forall\varphi\in C^\infty_0(\R)}$
is dense in $\dot H^{-\alpha}(0,L)$.

Consequently, (2.7) implies the left inequality (1.9).
The right inequality (1.9) is an immediate consequence of the observation
we made earlier: (1.7) and (1.8) imply that $R$ is a pseudodifferential
operator of order $-2\alpha=-(n+m)$.

Theorem 1.2 is proved.
\end{proof}

\begin{proof}[Proof of Theorem 1.3.]
Equations (1.4) and (1.7) imply
\begin{equation}
Ph=g:=QF.
\end{equation}
Here
\begin{equation}
F:=
\begin{cases} u_-, & x<0, \\ f, & 0\leq x\leq L,\\ u_+, & x>L, \end{cases}
\end{equation}
where
\begin{equation}
Qu_-=0, \quad x<0,
\end{equation}
\begin{equation}
Qu_+=0, \quad x>L,
\end{equation}
and $u_-$ and $u_+$ are chosen so that $F\in H^\alpha(\R)$.
This choice is equivalent to the conditions:
\begin{equation}
u^{(j)}_- (0) = f^{(j)}(0), \quad 0\leq j\leq \alpha-1,
\end{equation}
\begin{equation}
u^{(j)}_+(L)= f^{(j)}(L), \quad 0\leq j\leq \alpha-1.
\end{equation}
If $F\in H^\alpha(\R)$,
then $g:=QF\in H^{\alpha-n}(\R)=H^{-\frac{n+m}{2}}(\R)$,
and, by (2.9), one gets:
\begin{equation}
g=Qf+\sum^{n-\alpha-1}_{j=0}
\left[ a^-_j \delta^{(j)}(x) +a^+_j \delta^{(j)}(x-L) \right],
\end{equation}
where $a^\pm_j$ are some constants. There are
$n-\alpha=\frac{n+m}{2}$ constants $a^+_j$ and the same number of
constants $a^-_j$.

Let $G(x,y)$ be the fundamental solution of the equation
\begin{equation}
PG = \delta(x-y) \hbox{\ in\ } \R,
\end{equation}
which vanishes for $x<y$:
\begin{equation}
G(x,y)=0 \hbox{\ for \ }x<y.
\end{equation}

\begin{claim}
Such $G(x,y)$ exists and is unique. It solves the following Cauchy problem:
\begin{equation}
PG=0, \quad x>y, \quad G^{(j)}_x(x,y)\bigg|_{x=y+0}
=\delta_{j,m-1}, \quad 0\leq j\leq m-1,
\end{equation}
satisfies condition (2.16), and can be written as
\begin{equation}
G(x,y)=\sum^m_{j=1} c_j(y)\varphi_j(x), \quad x>y,
\end{equation}
where $\varphi_j(x)$, $1\leq j\leq m$, is a linearly independent
system of solutions to the equation:
\begin{equation}
P\varphi=0.
\end{equation}
\end{claim}
\begin{proof}[Proof of the claim.]
The coefficients $c_j(y)$ are defined by conditions (2.17):
\begin{equation}
\sum^m_{j=1} c_j(y) \varphi^{(k)}_j(y) = \delta_{k,m-1},
\quad 0\leq k\leq m-1.
\end{equation}
The determinant of linear system (2.20) is the Wronskian
$W(\varphi_1,\dots,\varphi_m)\not= 0$,
so that $c_j(y)$ are uniquely determined from (2.20).

The fact that the solution to (2.17), which satisfies (2.16), equals to 
the
solution to (2.15) -- (2.16) follows from the uniqueness of the
solution to (2.15) -- (2.16) and (2.17) -- (2.16),
and from the observation that the solution to (2.15) -- (2.16)
solves (2.17) -- (2.16).
The uniqueness of the solution to (2.17) -- (2.16) is a well-known result.

Let us prove uniqueness of the solution to (2.15) -- (2.16).
If there were two solutions, $G_1$ and $G_2$,
to (2.15) -- (2.16), then their difference $G:=G_1-G_2$,
would solve the problem:
\begin{equation}
PG=0 \hbox{\ in\ }\R,\quad G=0 \hbox{\ for\ }x<y.
\end{equation}
By the uniqueness of the solution to the Cauchy problem, it follows that
$G\equiv 0$. Note that this conclusion holds in the space of
distributions as well, because equation (2.21) has only the
classical solutions, as follows from the ellipticity
of $P$. Thus the claim is proved.
\end{proof}

From (2.8) and (2.14) -- (2.16) one gets:
\begin{align}
h=&\int^x_0 G(x,y) Qf\,dy
  + \int^x_0 G(x,y) \sum^{n-\alpha-1}_{j=0}
  \left[ a^-_j \delta^{(j)}(y) + a^+_j \delta^{(j)} (y-L) \right]dy 
\notag\\
=& \int^x_0 G(x,y) Qf\, dy + \sum^{n-\alpha-1}_{j=0} (-1)^j
  \left[ G^{(j)}_y (x,y)\bigg|_{y=0}
  a^-_j + G^{(j)}_y (x,y)\bigg|_{y=L} a^+_j \right]\notag \\
:=& \int^x_0 G(x,y) Qf\,dy+H(x).
\end{align}
It follows from (2.22) that $h\in H^{-\alpha}(\R)$
and
\begin{equation}
h=0 \hbox{\ for\ } x<0,
\end{equation}
that is, $(h,\varphi)=0$ $\forall \varphi\in C^\infty_0(\R)$ such that
$\supp \varphi \subset(-\infty,0)$.
In order to guarantee that $h\in\dot H^{-\alpha}(0,L)$ one has to satisfy
the condition
\begin{equation}
h=0 \hbox{\ for\ } x>L.
\end{equation}
Conditions (2.23) and (2.24) together are equivalent to
$\supp h \subset [0,L]$.
Note that although  $Qf\in \dot H^{-\frac{n+m}{2}}(0,L)$, 
so that $Qf$ is a distribution, the integral
$\int^x_0 G(x,y) Qf\,dy= \int^{\infty}_{-\infty} G(x,y) Qf\,dy$ 
is well defined as the unique solution
to the problem $Pw=Qf$, $w=0$ for $x<0$.

Let us prove that conditions (2.23) and (2.24) determine the constants
$a^\pm_j$, $0\leq j\leq \frac{n+m}{2}-1$, uniquely.

If this is proved, then Theorem 1.3 is proved, and formula (2.22) gives an
analytical solution to equation (1.4) in $\dot H^{-\alpha}(0,L)$
provided that an algorithm for finding $a^\pm_j$ is given.
Indeed, an algorithm for finding $G(x,y)$ consists of solving
(2.16) -- (2.17). Solving (2.16) -- (2.17) is accomplished
analytically by solving the linear algebraic system (2.20) and then
using formula (2.18). We assume that $m$ linearly
independent solutions $\varphi_j(x) $ to (2.19) are known.

Let us derive an algorithm for calculation of the constants
$a^\pm_j$, $0\leq j\leq \frac{n+m}{2}-1$,
from conditions (2.23) -- (2.24).

Because of (2.16), condition (2.23) is satisfied automatically by $h$ 
defined in (2.22).

To satisfy (2.24) it is necessary and sufficient to have
\begin{equation}
\int^L_0 G(x,y) Qf\,dy + H(x)\equiv 0 \hbox{\ for\ } x>L.
\end{equation}
By (2.18), and because the system $\{\varphi_j\}_{1\leq j \leq m}$ is
linearly independent, equation (2.25) is equivalent to the following set 
of equations:
\begin{equation}
\int^L_0 c_k(y) Qf\,dy +\sum^{\frac{n+m}{2}-1}_{j=0} (-1)^j
\left[ c^{(j)}_k(0) a^-_j + c^{(j)}_k(L) a^+_j \right]
=0, \quad 1\leq k\leq m.
\end{equation}
Let us check that there are exactly $m$ independent constants $a^\pm_j$
and that all the constants $a^\pm_j$ are uniquely determined by linear
system (2.26).

If there are $m$ independent constants $a^\pm_j$ and other constants
can be linearly represented through these, then linear algebraic system
(2.26) is uniquely solvable for these constants 
provided that the corresponding homogeneous system
has only the trivial solution. If $f=0$, then $h=0$, as follows from
Theorem 1.1,  and
$g=0$ in (2.14). Therefore $a^\pm_j=0$ $\forall j$,
and system (2.26) determines the constants $a^\pm_j$ $\forall j$
uniquely.

Finally, let us prove that there are exactly $m$ independent
constants $a^\pm_j$. Indeed, in formula (2.8) there are $\frac{n}{2}$
linearly independent solutions $u^-_j\in L^2(-\infty,0)$, so
\begin{equation}
u_-=\sum^{ \frac{n}{2}}_{j=1} b^-_j u^-_j,
\end{equation}
and, similarly, $u_+$ in (2.8) is of the form
\begin{equation}
u_+=\sum^{ \frac{n}{2}}_{j=1} b_j u^+_j,
\end{equation}
where $u^+_j\in L^2(0,\infty)$.
Condition $F\in H^\alpha(\R)$ implies
\begin{equation}
\sum^{ \frac{n}{2}}_{j=1} b^-_j(u^-_j)^{(k)} =f^{(k)} \hbox{\ at\ } x=0,
\quad 0\leq k\leq \alpha-1=\frac{n-m}{2}-1,
\end{equation}
and
\begin{equation}
\sum^{ \frac{n}{2}}_{j=1} b^+_j(u^+_j)^{(k)} =f^{(k)} \hbox{\ at\ } x=L,
\quad 0\leq k\leq \frac{n-m}{2}-1.
\end{equation}
Equations (2.29) and (2.30) imply that there are
$\frac{n}{2}-\frac{n-m}{2}=\frac{m}{2}$ independent constants
$b^-_j$ and $\frac{m}{2}$ independent constants $b^+_j$,
and the remaining $n-m$ constants $b^-_j$ and $b^+_j$ can
be represented through these $m$ constants by solving linear systems (2.29)
and (2.30) with respect to, say, first $\frac{n-m}{2}$ constants,
for example, for system (2.29), for the constants
$b^-_j$, $1\leq j\leq \frac{n-m}{2}$.
This can be done uniquely
because the matrices of the linear systems (2.29) and (2.30) are
nonsingular: they are  Wronskians of linearly independent solutions
$\{u^-_j\}_{1\leq j\leq \frac{n-m}{2}}$ and
$\{u^+_j\}_{1\leq j\leq \frac{n-m}{2}}.$ 

The constants 
$a_j^{\pm}$ can be expressed in terms of $b_j^{\pm}$ and $f$
by linear relations. Thus, there are exactly $m$ independent  
constants $a_j^{\pm}$.
This completes the proof of Theorem 1.3.
\end{proof}

\begin{remark}
In \cite[pp. 317-337]{R2} a theory of singular perturbations
for the equations of the form
\begin{equation}
\varepsilon h_\varepsilon + R h_\varepsilon = f
\end{equation}
is developed for a class of integral operators with a convolution
kernels $R(x,y)=R(x-y)$. This theory can
be generalized to the class of kernels $R(x,y)$ studied
in the present paper. The basic interesting problem is:
for any $\varepsilon>0$ equation (2.31) has a unique solution
$h_\varepsilon \in L^2(0,L)$; how can one
find the asymptotic behavior of $h_\varepsilon$ as $\varepsilon \to 0$?
The limit $h$ of $h_\varepsilon$ as $\varepsilon\to 0$ should
solve equation $Rh=f$ and, in general,
$h$ is a distribution, $h\in\dot H^{-\alpha}(0,L)$.
The theory presented in \cite[pp. 317-337]{R2} allows one
to solve the above problem for the class of kernels studied in the present 
paper.
\end{remark}
\begin{remark}
Theorems 1.1 and 1.2 and their proofs remain valid in the case when 
equation (1.4) is replaced by the
\begin{equation}
Rh:=\int_D R(x,y)h(y)dy=f,\quad x\in {\overline D}.
\end{equation}
Here $D\subset \R^r, \,r>1,$ is a bounded domain with a smooth boundary 
$S$, ${\overline D}$ is the closure of $D$,
$R(x,y)$ solves (1.7), where $P$ and $Q$ are uniformly elliptic
differential operators with smooth coefficients, $ord P=m\geq 0$,
$ord Q=n>m$, equation $Qh=0$ has only the trivial solution in 
$H^{\beta}(\R^r)$ for any fixed real number $\beta$.
Under the above assumptions, one can prove that the operator
defined by the kernel $R(x,y)$ is a pseudodifferential
elliptic operator of order $-2\alpha$, where $\alpha:=\frac {n-m}2$.
 We do not assume that $P$ and/or $Q$ are
selfadjoint or that $P$ and $Q$ commute. Therefore Theorems 1.1 and 1.2 
for equation (2.32) generalize the results from  \cite {R1}.
An analog of Remark 2.1 holds for the multidimensional equation (2.31) as 
well.
Equation (2.32) is the basic integral equation of random fields
estimation theory \cite {R1}.

In the short announcement by L.Piterbarg, Differentsialnye Uravneniya,
17, (1981), 2278-2279, (in Russian), some results concerning equation 
(2.32) are announced, but the 
author is not aware of the published proofs of these results.
\end{remark}

\section{Acknowledgement}
The author thanks Dr. K.Yagdjian for a discussion.

\end{document}